\documentclass{amsart}
\newtheorem{theorem}{Theorem}
\newtheorem{corollary}{Corollary}
\def\bibref[#1]{\cite{#1}}
\long\def\hh[#1,#2,#3]{h^{#1,#2}_#3}
\def\definition{\medskip\par \refstepcounter{equation}\noindent{\bf
Definition \theequation:} }
\def\remark{\medskip\par \refstepcounter{equation}\noindent{\bf Remark
\theequation:} }
\def\mod{\,\textup{mod}\,}
\def\braket#1{\langle#1\rangle}
\def\ket#1{|#1\rangle}

\def\C{\mathbb{C}}
\def\S{\mathbb{S}}
\def\N{\mathbb{N}}
\def\Z{\mathbb{Z}}
\def\calL{\mathcal{L}}
\def\P{\mathbb{P}}
\title{$n$-Schur Functions and Determinants on an Infinite Grassmannian}

\author{Alex Kasman}
\address{Mathematical Sciences Research Institute\\
Berkeley, CA 94720}
\begin{document}

\begin{abstract}
A set of functions is defined which is indexed by a positive
integer $n$ and  partitions of integers.  The case $n=1$ reproduces
the standard Schur polynomials.  These functions are seen to arise
naturally as a determinant of an action on the frame bundle of an
infinite grassmannian.  This fact is well known in the case of the
Schur polynomials ($n=1$) and has been used to decompose the
$\tau$-functions of the KP hierarchy as a sum.  In the same way, the
new functions introduced here ($n>1$) are used to expand
\textit{quotients} of $\tau$-functions as a sum with Pl\"ucker
coordinates as coefficients.
\end{abstract}

\maketitle

Among their many important properties, the Schur polynomials
\bibref[MacD] arise naturally as the determinant of an exponential
function acting on the frame bundle of the grassmannian of the Hilbert space
$H=L^2(S^1,\C)$ \bibref[PS].  It is for this reason that the
$\tau$-functions of the KP hierarchy can be expanded
as a sum of Schur polynomials with the Pl\"ucker coordinates as
coefficients \bibref[Sato,SW] (cf.\ \bibref[Engberg]).

\textit{Quotients} of $\tau$-functions have recently played a
prominent role in several papers on bispectrality \bibref[BHYbisp,KR],
Darboux transformations \bibref[BHYdt,aam] and random matrices
\bibref[AvM].  In \bibref[aam] these quotients themselves are
computed as a determinant of the action of a matrix valued function
on the frame bundle of the grassmannian $Gr^n=Gr(H^n)$ (cf.\
\bibref[PW]). 

This note will define the \textit{$n$-Schur functions} which play an
analogous role in this more general situation.  In particular, as in
the case $n=1$, it is shown that the determinant of the action of a
matrix function on the frame bundle of $Gr^n$ can be expanded as a sum
of $n$-Schur functions with the Pl\"ucker coordinates as
coefficients.  As an application one may expand quotients of
$\tau$-functions in this manner.

Note that although the $n$-Schur functions are defined in the next
section and only later shown to be related to determinants on the
frame bundle, this relationship should not be seen as a surprise.  It
is actually this property which led to the definition of these
functions.  The definition is given separately in the hope that
someone reading this paper, who might not have interest in infinite
grassmannians, may recognize other ways in which these generalized
Schur functions can be used.

\section{$n$-Schur Functions}

Let $n\in\N$ and consider the variables $\hh[i,j,k]$ ($1\leq i,j\leq
n$, $k\in\N$) which can be conveniently grouped into $n\times n$
matrices
$$
H_k=\left(\begin{matrix}
\hh[1,1,k]& \cdots& \hh[1,n,k]\cr
\vdots&\ddots&\vdots\cr
\hh[n,1,k]& \cdots& \hh[n,n,k]
\end{matrix}\right)\qquad k=0,1,2,3,4,\ldots
$$
Moreover, these matrices will be grouped into the infinite matrix
$M_{\infty}$ 
$$
M_{\infty}=\left(\begin{matrix}
\vdots&\vdots & \vdots& \cdots \cr
H_2&H_3&H_4&H_5&\cdots&\cdots\cr
H_1&H_2&H_3&H_4&\cdots&\cdots\cr
H_0&H_1&H_2&H_3&\cdots&\cdots  \hbox to 0pt{\hskip 1pc
$\leftarrow$\hbox{rows 0 through $n-1$}}\cr
0&H_0&H_1&H_2&H_3&\cdots\hbox to 0pt{\hskip 1pc
$\ \downarrow$\hbox{positively indexed rows}}\cr
0& 0&H_0&H_1&H_2&\cdots\cr
0& 0& 0&H_0&H_1&\cdots\cr
\vdots&&&&\ddots\end{matrix}\right)\qquad\qquad
$$
It is convenient to label the rows of this matrix by the integers with
$0$ being the first row with $H_0$ at the left and increasing downwards.

We wish to define a set of functions in these variables indexed by
partitions of integers.  Specifically, the index set $\S$ will be the
set of increasing sequences of integers whose values are eventually
equal to their indices
$$
\S=\{(s_0,s_1,s_2,\ldots)\ |\ s_{j+1}>s_j\in\Z\ \textup{and}\ \exists
N\ \textup{such\ that}\ j=s_j\ \forall j>N\}.
$$
Of particular interest here will be
the special case $0=(0,1,2,3,\ldots)\in\S$. 

\definition For any $S\in\S$ let $f_S^n:=\det(M_S\cdot M_0^{-1})$
where the infinite matrix $M_S$ is the matrix whose $j^{th}$ row is
the $s_j^{th}$ row of $M_{\infty}$. Then we have, for instance, that
$$
M_0=\left(\begin{matrix}
H_0&H_1&H_2&H_3&\cdots&\cdots\cr
0&H_0&H_1&H_2&H_3&\cdots\cr
0& 0&H_0&H_1&H_2&\cdots\cr
0& 0& 0&H_0&H_1&\cdots\cr
\vdots&&&&\ddots\end{matrix}\right).
$$
Equivalently, one could say that the element in 
position $(l,m)$ ($0\leq l,m\leq \infty$) of the matrix $M_S$ is given
by
$$
(M_S)_{l,m}=\hh[i,j,k]\qquad i=
1+(l\mod n),\ j=1+(s_m\mod n),\ c={\left\lfloor \frac
ln\right\rfloor-\left\lfloor \frac{s_m}n\right\rfloor}
$$
where $\hh[i,j,k]=0$ if $k<0$.

\remark Note that given any $N\in\N$ such that $s_i=i$ for all $i>Nn$,
 the matrix $M_S\cdot M_0^{-1}$ looks like the identity matrix below
 the $Nn^{th}$ row.  Consequently, the easiest way to actually compute
 these functions is as two finite determinants
$$
f_S^n=\frac{\det(M_S|_{nN\times nN})}{(\det H_0)^N}
$$
where $M_S|_{nN\times nN}$ denotes the top left block of size
$nN\times nN$ of the matrix $M_S$.  For example,
regardless of $n$, $f_0^n\equiv 1$.  However, for
$S=(-2,1,2,3,\ldots)$ one finds instead 
$$
f_S^1=\frac{\hh[1,1,2]}{\hh[1,1,0]}\qquad
f_S^2=\frac{\hh[1,1,1]\hh[2,2,0]-\hh[1,2,0]\hh[2,1,1]}{\hh[1,1,0]\hh[2,2,0]-\hh[1,2,0]\hh[2,1,0]}.
$$

\remark In fact, in the case $n=1$ and
$\hh[1,1,0]=1$, the functions $\{f_S^1\}$ are the famous Schur polynomials
\bibref[MacD] (cf.\ \bibref[PS,SW]).
Similarly, if we assume in general that
$\det H_0=1$, then
all $f_S^n$ are polynomials.   

\remark If we consider $\hh[i,j,k]$ to have weight $kn+i-j$ then the
function $f_S^n$ is homogenous of weight $\sum_{j=0}^{\infty}s_j-j$.

\section{The Grassmannian $Gr^n$}

In this section we will recall notation and some basic facts about an
infinite dimensional grassmannian.  Please refer to
\bibref[aam,PS,PW,SW] for additional details.

Let $H^n=L^2(S^1,\C)$ be the Hilbert space of square-integrable vector
valued functions $S^1\to\C^n$, where $S^1\subset\C$ is the unit
circle.  Denote by $e_i$ ($0\leq i \leq n-1$) the $n$-vector which has
the value 1 in the $i+1^{st}$ component and zero in the others.  We
fix as a basis for $H^n$ the set $\{e_i|i\in\Z\}$ with
$$
e_i:=z^{\lfloor \frac i n \rfloor}e_{(i\mod n)}.
$$
The Hilbert space has the decomposition
\begin{equation}
H^n=H_+^n\oplus H_-^n\label{eqn:split}
\end{equation}
where these subspaces are spanned by the basis elements with
non-negative and negative indices respectively.
Then $Gr^n$ denotes the grassmannian of all closed subspaces $W\subset
H^n$ such that the orthogonal projection $W\to H^n_-$ is a compact
operator and such that the orthogonal projection $W\to H_+^n$ is
Fredholm of index zero \bibref[PS,PW].

Associate to any basis $\{w_0,w_1,\ldots\}$ for a point $W\in Gr^n$ 
the linear map $w$
\begin{eqnarray*}
w:H_+^n&\to&W\\
e_i&\mapsto&w_i.
\end{eqnarray*}
The basis is said to be \textit{admissible} if
$w$ differs from the identity by an element of trace class
\bibref[Simon].  The \textit{frame bundle} of $Gr^n$ is the set of
pairs $(W,w)$ where $W\in Gr^n$ and $w:H_+^n\to W$ is an admissible
basis.

There is a convenient way to embed $Gr^n$ in a projective space.
Let $\Lambda$ denote the infinite alternating exterior algebra
generated by the alternating tensors
$$
\{e_{s_0}\wedge e_{s_1}\wedge
e_{s_2}\wedge\cdots|(s_0,s_1,s_2,\ldots)\in\S\}.
$$
To any point $(W,w)$ in the frame bundle we associate the alternating
tensor
$$
\ket{w}:=w_0\wedge w_1\wedge w_2\wedge\cdots\in\Lambda.
$$
Note in particular that $\ket{\cdot}$ is projectively well defined on
the entire fiber of $W$ (i.e.\ for two admissible bases of $W$ we have
$\ket{w}=\lambda\ket{w'}$ for some non-zero
constant $\lambda$).  Consequently, $\ket{W}$ is a well defined
element of the projective space $\P\Lambda$.  

The Pl\"ucker coordinates (cf.\ \bibref[HodgePedoe]) of $W$ are
defined as the coefficients $\braket{S|W}$ in the unique expansion
$$
\ket{W}=\sum_{S\in\S} \braket{S|W} e_{s_0}\wedge e_{s_1}\wedge
e_{s_2}\wedge \cdots
$$
and are therefore well defined as a set up to a common multiple.
Alternatively, given an admissible basis $w$ for $W$, $\braket{S|W}$
is the determinant of the infinite matrix made of the rows of $w$
indexed by the elements of $S$. 

\section{Determinants of Action on Frame Bundle}

Let $g$ be an $n\times n$ matrix valued function of $z$ with expansion
\begin{equation}
g=\sum_{k=0}^{\infty} H_k z^k\qquad H_k=\left(\begin{matrix}
\hh[1,1,k]& \cdots& \hh[1,n,k]\cr
\vdots&\ddots&\vdots\cr
\hh[n,1,k]& \cdots& \hh[n,n,k]
\end{matrix}\right)\label{eqn:gform}
\end{equation}
such that an inverse matrix $g^{-1}$ exists for all $z$.  We will view
$g\in GL(H^n)$ as an operator on $Gr^n$ and demonstrate that the
$n$-Schur functions arise naturally in this context.

In general, an operator on the frame bundle \bibref[PS] is a pair
$A=(g,q)$ where $g\in GL(H^n)$ with the form 
\begin{equation}
g=\left(\begin{matrix}a&b\cr
c&d\end{matrix}\right)
\end{equation}
relative to the splitting \eqref{eqn:split} and $q:H_+^n\to H_+^n$
such that $a\cdot q^{-1}$ differs from the identity by an operator of
trace class.  The action is given by
$$
A:(W,w)\mapsto (gW,gwq^{-1}).
$$

In the particular case \eqref{eqn:gform} of interest here 
$c=0$ and we simply let $q=a$ so that $aq^{-1}$ is
the identity matrix.  We will write $\ket{g|w}=\ket{gwa^{-1}}$ for the
action of $g$ on the frame bundle.  Moreover, since this action is
well defined on projective equivalence classes we will write
$\ket{g|W}$ for the class containing $\ket{g|w}$ with any admissible
basis $w$ of $W$.

The main result of this paper is then the observation that for any
point $W\in Gr^n$ the determinant $\braket{0|g|W}$ can be written in
terms of the Pl\"ucker coordinates of $W$ and the $n$-Schur functions:
\begin{theorem}\label{thm}
For $W\in Gr^n$ and $g$ as in \eqref{eqn:gform}
$$
\braket{0|g|W}=\sum_{S\in\S}\braket{S|W} f_S^n.
$$
\end{theorem}

The proof is elementary in the case $W=W_S$ with basis
$\{e_{s_0},e_{s_1},e_{s_2},\ldots\}$.  In fact, this is essentially
the definition of $f_S^n$ since the matrix representation of the
operator $g$ is precisely $M_{\infty}$.  The general case follows from
the observation that multilinearity of determinants is equivalent to
the linearity of the map $\langle 0|g|:\Lambda\to\C$ and expanding
$\ket{W}$ as a sum.

\section{Application: Expanding Quotients of $\tau$-functions}

\subsection{Basic Facts about the KP Hierarchy}
A pseudo-differential operator of the form
\begin{equation}
\calL=\partial + u_1(x)\partial^{-1}+ u_2(x)\partial^{-2}+\cdots\qquad
\partial=\partial/\partial x\label{eqn:initcond}
\end{equation}
is said to be a solution of the KP hierarchy if its coefficients $u_i$
depend on ``time variables'' $t_1,t_2,\ldots$ so as to satisfy the equations
\begin{equation}
\frac{\partial}{\partial t_i}\calL=[\calL,(\calL^i)_+]\label{eqn:KPhier}
\end{equation}
where the ``$+$'' subscript indicates projection onto the differential
operators by simply eliminating all negative powers of $\partial$ and
$[A,B]=A\circ B- B\circ A$ \bibref[SW]. 
Remarkably, there exists a convenient way to encode all information
about the KP solution $\calL$ in a single function of the time
variables $t_1,t_2,\ldots$.  Specifically, each of the coefficients
$u_i$ of $\calL$ can be written as a rational function of this
function $\tau(t_1,t_2,\ldots)$ and its derivatives \bibref[SW].
Alternatively, one can construct $\calL$ from $\tau$ by letting $W$ be
the pseudo-differential operator
$$
W=\frac{1}{\tau}
\tau(t_1-\partial^{-1},t_2-\frac12\partial^{-2},\ldots)
$$
and then $\calL:=W\circ\partial\circ W^{-1}$ is a solution to the KP
hierarchy.  
 
\subsection{Expanding $\tau$-quotients with $n$-Schur functions}
Of particular interest for many situations are the solutions $\calL$
of the KP hierarchy that have the property that
$L=\calL^n=(\calL^n)_+$ is an \textit{ordinary} differential operator.
We say that solutions of the KP hierarchy with this property are
solutions of the of the $n$-KdV hierarchy.
Associated to any chosen solution $L=\calL^n$ of the $n$-KdV hierarchy, we
naturally associate an $n$-vector valued function of the variables
$\{t_i\}$ and the new \textit{spectral parameter} $z$.  In particular,
following \bibref[PW] we define the \textit{vector Baker-Akhiezer}
function to be the unique function $\vec\psi(z,t_1,t_2,\ldots)$
satisfying 
\begin{equation}
\calL^n\vec\psi=L\vec\psi=z\vec\psi\qquad \frac{\partial}{\partial
t_i}\vec\psi=(\calL^i)_+ \vec\psi\label{eqn:psivec}
\end{equation}
and such that the $n\times n$
Wronskian matrix 
\begin{equation}
\Psi(z,t_1,t_2,\ldots):=\left(\begin{matrix}
\vec\psi\cr
\frac{\partial}{\partial x} \vec\psi
\cr
\vdots\cr
\frac{\partial^{n-1}}{\partial x^{n-1}} \vec\psi\end{matrix}\right)\label{eqn:wr}
\end{equation}
is the identity matrix when evaluated at
$0=t_1=t_2=\ldots$.

Defining $g=\Psi^{-1}\in GL(H^n)$ where $\Psi$ is the matrix above,
the determinant $\braket{0|g|W}$ is a (projective) function of the
variables $t_i$.  
It is shown in \bibref[aam] (Definition 7.4 and Claim 7.12) that these functions are  quotients of KP
$\tau$-functions with a $\tau$-function corresponding to $\calL$ in
the denominator.  Consequently, we may write these quotients in terms
of the $n$-Schur functions:
\begin{corollary}
Let $\calL$ be a solution of the $n$-KdV hierarchy with corresponding
 matrix $\Psi$ given in \eqref{eqn:psivec} and \eqref{eqn:wr}.
Defining time dependent variables $\hh[i,j,k]$ by
$$
\Psi^{-1}=\sum_{k=0}^{\infty} H_k z^k.
$$
gives the $n$-Schur functions dependence on the time variables of the
KP hierarchy.
Then one has that
every $f_S^n$ is then a quotient of KP $\tau$-functions in the sense
of \bibref[aam] for any $S\in\S$.  Moreover, it follows that
$$
\tau_0\cdot(\sum_{S\in\S} \pi_S f_S^n)
$$
is a $\tau$ function of the KP hierarchy whenever $\tau_0$ is a
$\tau$-function for $\calL$ and $\pi_S$ are the Pl\"ucker coordinates
of some point in $Gr^n_{rat}$.
\end{corollary}

The special case that $\tau_0=1$, $\calL=\partial$ and $\Psi=\exp(\sum
t_iz^i)$ reproduces the well known expansion of KP $\tau$-functions in
terms of the Schur polynomials \bibref[SW] (see also
\bibref[Engberg]).  Explicit new examples based on the results
presented here, including ``higher rank'' \bibref[PW] cases, can be
found in \bibref[Thesis] and \bibref[GeominDyn] where the $n$-Schur
functions are used to demonstrate the connection between the KP
hierarchy and Pl\"ucker relations \bibref[HodgePedoe].

\medskip
\par
\noindent{Acknowledgements:} The new results reported here were mostly
obtained as part of my Ph.D. thesis \bibref[Thesis].  I am extremely
grateful to my advisor, E. Previato, for all of her assistance and
guidance.  I also appreciate many helpful discussions with F. Sottile and
M. Bergvelt.

\end{document}